\documentclass[10pt]{article}
\usepackage{hyperref}
\usepackage{amsmath}
\usepackage{amsthm}
\usepackage{epsfig}
\usepackage{amssymb}
\usepackage{dsfont}
\usepackage{graphicx}
\usepackage{float}

\newtheorem{theorem}{Theorem}[section]
\newtheorem{lemma}[theorem]{Lemma}
\newtheorem{mthm}{Theorem}

\newtheorem{de}{Definition}[section]

\catcode`@=11
\@addtoreset{equation}{section}
\catcode`@=12

\begin{document}

\title{A degree theory approach for the shooting method}
\author{Congming Li}
\date{}
\maketitle

\begin{abstract}

The classical shooting-method is about finding a suitable initial shooting positions
to hit the desired target.
The new approach formulated here, with the introduction and  the analysis of
the `target map' at its core,
naturally connects the classical shooting-method
to the simple and beautiful topological degree theory.

We apply the new approach, to a motivating example,
to derive the existence of global positive solutions of
the Hardy-Littlewood-Sobolev (also known as Lane-Emden) type system:
\[
\left\{ \begin{aligned}
        &(-\triangle)^ku(x) = v^p(x),  \,\, u(x)>0 \quad\text{in}\quad\mathbb{R}^n,\\
         &         (-\triangle)^k v(x) =u^q(x), \,\, v(x)>0 \quad\text{in}\quad\mathbb{R}^n,  p, q>0,
         \end{aligned} \right.
\]
in the critical and supercritical cases $\frac{1}{p+1}+\frac{1}{q+1}\leq\frac{n-2k}{n}$.
Here we derive the existence with the computation of the topological degree
of a suitably defined target map. This and some other results presented in this article completely solve some
long-standing open problems about the existence or non-existence of positive entire solutions.

\end{abstract}



\section{Introduction}

The degree theory approach for the shooting method consists of three
components: (a) defining the target map with the careful choice of the domain of possible initial shooting positions and the range
of all possible targets;
(b) analyzing  the target map and showing that the map is onto or has fixed points via the degree theory;
(c) proving that the fixed point or the special target obtained in (b) leads to
solutions of our partial differential equations or dynamical systems.

We apply the new approach to show the existence of global positive solutions to:
\begin{equation} \label{eq:1}
\left\{ \begin{aligned}
         &-\triangle u_i = f_i(u) \quad\text{in}\quad\mathbb{R}^n,\quad i=1,\ldots,L\\
          &      u_i>0\quad\text{in}\quad\mathbb{R}^n. 
                          \end{aligned} \right.
\end{equation}

We rewrite the above system in radial coordinate as an initial value problem and choose
a set of suitable initial values as the domain of our target map and then give a suitable definition of the target map.
With the proof of the continuity of the target map, we
apply the  degree theory as presented in \cite{Nirenberg01} and \cite{Schwartz69}
 to compute the index of the target map and
to show that the target map is onto. We then prove that this guarantees
the existence of some global positive solutions.

We only need some very mild assumptions on the function: $F=(f_1,\ldots,f_L):\mathbb{R}_+^L=\underbrace{\mathbb{R}_+
\times\mathbb{R}_+\ldots\times\mathbb{R}_+}_L\rightarrow \mathbb{R}_+^L$, where $\mathbb{R}_+=[0,\infty)$.
In this article,
we \textbf{always} assume that $F$ is continuous on $\mathbb{R}_+^L$ and is
Lipschitz in the interior of $\mathbb{R}_+^L$.

\begin{de}
We say that $F$, or system (\ref{eq:1}), is \emph{non-degenerate} or  \emph{non-reducible}
 if for any permutation $i_1,\ldots,i_L$ of
$1,2,\ldots,L$, any $1 \leq k <L$, and any $u\in \mathbb{R}_+^L$ with $u_{i_1}>0,\ldots,u_{i_k}>0$ and $u_{i_{k+1}}=\ldots=u_{i_L}=0$, we have
$f_{i_{k+1}}(u)+\ldots+f_{i_L}(u)>0$. 
\end{de}

A trivial degenerate case is the following:

\[ \left\{ \begin{array}{l}
         -\triangle u_1 = u_1^p, \\
                  -\triangle u_2=u_2^p.
                  \end{array} \right.
         \]
The above system can be decoupled to the study of two scalar equations.

The main idea is to associate the existence of solutions to (\ref{eq:1}) with the non-existence
of solutions to the Dirichlet boundary value problem of the same elliptic system on balls:
\begin{equation} \label{eq:4}
\left\{ \begin{aligned}
         &-\triangle u_i = f_i(u) \quad\text{in}\quad B_R\subset\mathbb{R}^n, i=1,\ldots,L\\
                 & u_i>0  \quad\text{in} \quad B_R\\
                 & u_i=0 \quad\text{on} \quad\partial B_R,
                          \end{aligned} \right.
                          \end{equation}
$B_R=B_R(0)=\{x\in\mathbb{R}^n||x|<R \}$.

\begin{mthm}\label{thm:1}
Assuming that $F$ is non-degenerate, then system (\ref{eq:1}) admits a solution if
the corresponding system (\ref{eq:4}) admits no solution for any given $R>0$.
Furthermore, if we assume that $F(u)\neq 0$
for $u>0$  (we write $u>0$ if all components of $u$ are positive i.e. $u_i>0$ for $i=1,\ldots,L$),
then (\ref{eq:1}) admits a solution $u(x)$ such that
$u(x) \rightarrow 0$ uniformly as $|x| \rightarrow \infty $.
\end{mthm}

The above theorem is a consequence of the following theorem:
\begin{mthm}\label{thm:2}
Assuming that $F$ is non-degenerate, then system (\ref{eq:1})
admits a radially symmetric solution if system (\ref{eq:4})
admits no radially symmetric solutions.
Furthermore, if we assume that $F(u)\neq 0$ for $u>0$ and small,
then (\ref{eq:1}) admits a radially symmetric solution $u(r)$ such that
$u(r) \rightarrow 0$ uniformly as $r \rightarrow \infty $.
\end{mthm}

A motivating example as well as an important application of Theorem \ref{thm:1} is:
\begin{mthm}\label{thm:3}
The system:
\begin{equation} \label{eq:6}
\left\{ \begin{aligned}
       &(-\triangle)^ku = v^p, \,\, p>0 \quad\text{in}\quad\mathbb{R}^n,\\
         &(-\triangle)^k v =u^q, \,\, q>0 \quad\text{in}\quad\mathbb{R}^n,\\
                  &u,v>0,\quad\text{in}\quad\mathbb{R}^n,\\
                  &u(x), \, v(x) \rightarrow 0 \text{ uniformly as } |x| \rightarrow \infty,
                          \end{aligned} \right.
\end{equation}
admits a positive solution in the critical and super-critical cases
$\frac{1}{p+1}+\frac{1}{q+1}\leq\frac{n-2k}{n}$
for any $2k<n$.
\end{mthm}
This system can be regarded as the `blow up' equations for a large class of systems of
nonlinear equations arising from geometric analysis
and other physical sciences.

The Lane-Emden conjecture (see \cite{CL09}, \cite{CLO06}, \cite{LL13preprint}, \cite{SZ96}-\cite{Souplet09}),
which was proved in many cases, states that:
in the subcritical cases $\frac{1}{p+1}+\frac{1}{q+1}> \frac{n-2k}{n}$, system (\ref{eq:6})
admits no solutions. The complete resolution of the Lane-Emden conjecture
and our result here provide a complete understanding of existence and non-existence of global positive solutions of (\ref{eq:6}).

For the corresponding scalar case with
$u=v$ and $p=q$, we have the following result:

\begin{mthm} \label{thm:4}
For $n>2k$, the problem:
 \begin{equation} \label{eq:scalar}
\left\{ \begin{aligned}
        &(-\triangle)^ku = u^p \\
         &         u>0, \quad\text{ in }\quad\mathbb{R}^n
                          \end{aligned} \right.
\end{equation}
admits a solution if and only if $p \geq \frac{n+2k}{n-2k}$. It admits a finite energy solution
in the sense $\int_{\mathbb{R}^n}u^{p+1}<\infty$ if and only if $p = \frac{n+2k}{n-2k}$.
\end{mthm}
This complete resolution of the exsitence/nonexistence is very surprising and is
 the best indication of the power of the degree theory approach for the shooting method.

For equations (\ref{eq:1})-(\ref{eq:scalar}), we introduce  the target map
and then apply the degree theory (property 1.5.5 in \cite{Nirenberg01}) to
show that the target map is onto which guarantees that we can shoot to the desired target.
The existence of global positive solutions follows from this.

We start with some background materials.
The existence of (\ref{eq:scalar})
in the critical and super-critical cases $p \geq \frac{n+2k}{n-2k}$ and
the nonexistence in the subcritical case have
been investigated extensively in the last 30 years (see \cite{BrezisNirenberg83}-\cite{CL97},
\cite{Dancer01}, \cite{Gui96}, \cite{NS86}, \cite{PS86}, \cite{PV08}, and \cite{SZ96}-\cite{Souplet09}).
When $k=1$, one can use the classical shooting method (\cite{CEQ94} and \cite{JL73})
to show the existence of solutions to:
\begin{equation} \label{eq:scalar1.5}
\left\{ \begin{aligned}
        &-\triangle u = u^p \\
         &         u>0,  \quad\text{ in }\quad\mathbb{R}^n
                          \end{aligned} \right.
\end{equation}
in the critical and super-critical cases $p\geq \frac{n+2}{n-2}$ and $n \geq 3$.

In fact, one seeks radially symmetric solutions
$u(x)=w(|x|)$
and writes the above equation as a second order ordinary differential equation  in $w$
with initial value $w(0)=\alpha >0$ and $w'(0)=0$. In the critical and super-critical case,
with the non-existence of solution of the same equation on a ball with Dirichlet boundary value,
one sees that $w(r)>0$
for all $r>0$ and $w (r) \searrow 0$ as $r \longrightarrow \infty$. Basically,
the solutions corresponding to different initial values are just a scaling
change of each other.

When $k \geq 2$, instead of one dimensional initial value which scales to each other,
we are encountered with multi-dimensional initial value:
$$\alpha=(\alpha_1, \alpha_2, \cdots, \alpha_{k-1}) \in R^{k-1}$$ with
$ \alpha_i=(-\triangle)^{i-1} u(0)$, for $i=1,2,...k-1$. Among them,
in the critical cases as well as in many super-critical cases,
there is at most one scaling class (one-dimensional) of initial values
from which we can shoot to a global solution. To show the existence of positive solutions of
(\ref{eq:scalar}), up to a simple scaling, we have to find the special 1-D initial values.
This is the main reason why there are so many results  about (\ref{eq:scalar1.5})
but very little about (\ref{eq:scalar}) for a long time period.

One  notices that (\ref{eq:scalar}) has a variational structure.
It is then very natural to seek a critical point and thus find
a solution. This is indeed what has been done in the critical case
$\frac{1}{p+1}+\frac{1}{q+1} = \frac{n-2k}{n}$
in \cite{Lieb83}.
Unfortunately,  in both the super-critical
and the subcritical cases, Lei and Li showed in \cite{LL13preprint} that
system (\ref{eq:scalar}) or system (\ref{eq:6}) has no solution with finite energy.

The degree theory approach for the shooting method gives a simple solution to this difficult problem.
It can be used to solve a much large class of problems. We first derive theorems \ref{thm:1} and \ref{thm:2}
for the general but simple looking
system (\ref{eq:1}). We then apply the general results in theorems \ref{thm:1} and \ref{thm:2}
to the existence
of global positive solutions for the Hardy-Littlewood-Sobolev type systems (\ref{eq:6}) and (\ref{eq:scalar})
in theorems \ref{thm:3} and \ref{thm:4}.

First, if $F(\beta)=0$ for some $\beta>0$, then
$u=\beta$ is a trivial/stationary/constant/ positive solution of (\ref{eq:1}) and
we are done with our proof of theorem \ref{thm:1}.

Thus, we may assume that $F(\beta)\neq 0$ for $\beta>0$.
Later, we can see that we only need the assumption holds for $\beta$ small.
The key ingredient of this paper is the definition and the analysis of the
`target map' $\psi$.

For any initial value $\alpha=(\alpha_1, \cdots, \alpha_L)$ with $\alpha_i>0, i=1,2,....L$,
we solve the following initial value
problem and denote the solution as $u(r,\alpha)$:
\begin{equation} \label{eq:ode}
\left\{ \begin{aligned}
        &u_i^{''}(r) +\frac{n-1}{r} u_i^{'}(r)= -f_i(u) \\
         &u_i'(0)=0, u_i(0)=\alpha_i \quad i=1,2,\ldots,L.
                          \end{aligned} \right.
\end{equation}

For $\alpha > 0$, we define the target map with $\psi(\alpha)=u(r_0, \alpha)$
where $r_0$ is the smallest value of $r$ for which
$u_i(r, \alpha)=0$ for some $i$ or when there exists no such $r$, we define
$\psi(\alpha)=\lim_{r \longrightarrow \infty} u(r, \alpha)$. In the later case, one can
see that $F(\psi(\alpha))=0$.
This and the assumption that $F(\beta) \neq 0$ when $\beta>0$
ensure that $\psi(\alpha) \in \partial R_+^L$.
With the natural extension of $\psi(\alpha)=\alpha$ for
$\alpha \in \partial R_+^L$, we then show that $\psi$ is continuous from $R_+^L$
to $ \partial R_+^L$.

Applying the degree
theory, we show that $\psi$ is onto from $A_a$ to $B_a$ where:

\begin{equation} \label{eq:scalar1}
\left\{ \begin{aligned}
 & A_a\triangleq \{ \alpha \in R_+^L \,\, {\bf \mid } \quad {\tiny \displaystyle\sum_{i=1,\cdots, L}} \alpha_i=a\},\\
  &B_a\triangleq\{ \alpha \in \partial R_+^L \,\, {\bf \mid } \quad {\tiny \displaystyle\sum_{i=1,\cdots, L}} \alpha_i \leq a\},
                          \end{aligned} \right.
\end{equation}
for any $a>0$.
In particular, there exists at least one
$\alpha_a \in A_a$ for
every $a>0$ such that $\psi(\alpha_a)=0$.

Shooting from the initial value $\alpha_a$, using the fact that the system
(\ref{eq:4}) admits no radially symmetric solution, we obtain a solution of (\ref{eq:1}).
In fact, we get a solution for every $a>0$.
We remark that we get infinity many solutions
even if our assumptions on $F$ only hold for $u$ small.
In this case, we just employ the above method for $a$ small.

Obviously, our method can also be employed to study more general dynamic systems of the form:
\[
\left\{ \begin{aligned}
         \frac{dU}{dt}&=F(t,U), \\
                  U(0)&=\alpha \in \mathbb{R}_+^L,
                          \end{aligned} \right.
                          \]
and to seek a initial state $U(0)=\alpha$ from which we get a global non-negative solution:
$U(t): [0,\infty) \longrightarrow \mathbb{R}_+^L$.

\section{Proofs of theorems \ref{thm:1}-\ref{thm:4}}

It is clear that theorem \ref{thm:1} is a consequence of theorem \ref{thm:2}.\\

\textbf{Proof of theorem \ref{thm:2}}: \\

As we have discussed in the previous section,  we assume $F(\beta) \neq 0$ when $\beta > 0$
and define $\psi(\alpha)$ to be the target map from $\mathbb{R}_+^L$ to $\partial \mathbb{R}_+^L$.

\begin{lemma} \label{lemma:1}
The map $\psi$:
$\mathbb{R}_+^L \rightarrow \partial \mathbb{R}_+^L$ is continuous.
\end{lemma}

We prove this lemma later. Instead, we apply the degree theory to show that:

\begin{lemma} \label{lemma:2}
For any $a>0$, $\psi$ is an onto map from $A_a$ to $B_a$ and  thus
there exists at least one $\alpha_a \in A_a$
such that $\psi(\alpha_a)=0$.
\end{lemma}

{\bf Proof of lemma \ref{lemma:2}}\\

Recall that $B_a= \{ \alpha \in \partial \mathbb{R}_+^L \quad {\bf |} \quad
\displaystyle\displaystyle\sum_{i=1,\cdots, L} \alpha_i \leq a\}$, then as a consequence of the non-increasing property
of the solutions of (\ref{eq:ode}), we see that $\psi$ maps $ A_a \longrightarrow B_a$.

Let $\phi(\alpha)=\alpha+\frac{1}{L}(a- \displaystyle\sum_{i=1,\cdots, L} \alpha_i)(1,\cdots,1): B_a \longrightarrow A_a$,
then $\phi$ is continuous with a continuous inverse $\phi^{-1}(\alpha)=\alpha-(\displaystyle\min_{i=1,\cdots, L} \alpha_i)(1,\cdots,1): A_a \longrightarrow B_a$.

The map: $G=\phi \circ \psi: A_a \longrightarrow A_a$ is continuous and $G(\alpha)=\alpha$
on the boundary of $A_a$. Based on the Heinz-Lax-Nirenberg version of the degree theory,
property 1.5.5 in page 8 of \cite{Nirenberg01},
we calculate that $deg(G, A_a, \alpha)=deg(Identity, A_a, \alpha)=1$ for any interior point
$\alpha \in A_a$. Consequently, $G$ is onto which implies that $\psi$ is also onto. This shows that there
exists an $\alpha_a \in A_a$
such that $\psi(\alpha_a)=0$ for any $a>0$.
This completes the proof of lemma \ref{lemma:2}.\\

{\bf Showing the existence of solutions:} \\

First, we show that the solution $\overline{u}(r)$ of (\ref{eq:ode}) with the initial value
$\overline{u}(0)=\alpha_a$ never touches the wall for finite $r$ and thus is defined for all $r>0$.
Suppose in the contrary that
$r_0$ is the first positive number such that $\overline{u_i}(r_0)=0$ for some $i$. Then by definition
$\overline{u}(r_0)=\psi(\alpha_a)=0$. This implies that
$u(x)=\overline{u}(|x|)$ is a solution of
(\ref{eq:4}) with $R=r_0$.
This contradicts with the assumption that system (\ref{eq:4}) admits no solutions.
Consequently, we get $\overline{u}_i(r)>0$ for $i=1,\cdots, L$ and $r>0$
and $\displaystyle\lim_{r \rightarrow \infty}\overline{u}(r)=\psi(\alpha_a)= 0$.

Clearly, $u(x)=\overline{u}(|x|)$ is an radially symmetric classical solution of
(\ref{eq:1}) with $u(x) \rightarrow 0$ uniformly as $|x| \rightarrow \infty $.

Thus to complete the proof of theorem \ref{thm:2}, we only need to establish lemma \ref{lemma:1}:
{\it for any $\overline{\alpha} \in \mathbb{R}_+^L$,
$\psi$ is continuous at $\overline{\alpha}$}.\\

{\bf Proof of lemma \ref{lemma:1}:} \\

There are three cases to be considered:
\begin{enumerate}
\item $\overline{\alpha} \in \partial \mathbb{R}_+^L$.

\item  $\overline{\alpha}>0$,  and the solution $u(r,\overline{\alpha})$ of
(\ref{eq:ode}) with initial value $\overline{\alpha}$
touches the wall at the smallest possible value $r_0$ with $u_{i_0}(r_0, \overline{\alpha})=0$,
for some $1\leq i_0\leq L$.

\item $\overline{\alpha}>0$,
and the solution $u(r, \overline{\alpha})$ of
(\ref{eq:ode}) never touches the wall, or $u_i(r,\overline{\alpha})>0$ for $i=1,\ldots, L$ and $ r \in [0,\infty)$.

\end{enumerate}

\textbf{Case (1)}:\\

If $\overline{\alpha}=0$, then $|\psi(\alpha) - \psi(\overline{\alpha})|= |\psi(\alpha)|
\leq |\alpha | = |\alpha-\overline{\alpha}|
\rightarrow 0$ as $\alpha \longrightarrow 0=\overline{\alpha}$.

When $\overline{\alpha}\neq0$, without loss of generality, we assume that $\overline{\alpha_1}=0, \ldots, \overline{\alpha_j}=0$
and $\overline{\alpha_{j+1}}>0, \ldots, \overline{\alpha_L}>0$, $1\leq j <L$. We must have $f_1(\overline{\alpha})+\cdots+
f_j(\overline{\alpha}) >0$ by the assumption that $F$ is non-degenerate. Thus we may assume $f_1(\overline{\alpha})=c>0$.
By continuity, there exists $\delta_1>0$ such that
$f_1(\alpha)\geq \frac{c}{2}$ if $|\alpha-\overline{\alpha}| \leq \delta_1$. Classical ODE theory
shows that there exists a $\delta_2>0$ such that if $|\alpha-\overline{\alpha}| \leq \delta_2$
and $r<\delta_2$
then $|u(r,\alpha)-\overline{\alpha}| \leq \delta_1$ before $u(r,\alpha)$ touches the wall.
One sees that as $|\alpha-\overline{\alpha}| \longrightarrow 0$, $u_1(r_{\alpha},\alpha)=0$
for some $r_{\alpha} \longrightarrow 0$. Hence,
$|\psi(\alpha) - \psi(\overline{\alpha})| \leq |u(r_{\alpha}, {\alpha})-\alpha| +|\alpha-\overline{\alpha}|
\rightarrow 0$ as $\alpha \longrightarrow\overline{\alpha}$.\\

\textbf{Case (2)}:\\

From the fact that $f_i \geq 0$, one derives that $ u^{'}_{i_0}(r_0, \overline{\alpha})<0$. This transversality
condition and the ODE stability imply that $\psi$ is continuous at $\overline{\alpha}$.\\

\textbf{Case (3)}:\\

In this case, we first show that $\psi(\overline{\alpha})=0$.
Classical ODE or Elliptic theory shows that $F(\psi(\overline{\alpha}))=0$. Hence $\psi(\overline{\alpha})
\in \partial \mathbb{R}_+^L$.
By our non-degeneracy assumption, we conclude that $\psi(\overline{\alpha})=0$. In fact, the
 non-degeneracy condition implies that $F(\beta) \neq 0$ when $\beta \in \partial \mathbb{R}_+^L$ and
$\beta  \neq 0$.

Then $u(r,\overline{\alpha})$ is positive and small for $r$ large. Continuous dependence of initial values for our ODE
implies that for any $R$ large but fixed when $\alpha$ is close to $\overline{\alpha}$, then
$u(r,\alpha)>0$ for $r \in [0,R]$ and $u(R,\alpha)$ is close to $u(r,\overline{\alpha})$ and thus is small.
Consequently $|\psi(\alpha)| \leq |u(R,\alpha)|$ is small. This shows that $\psi$ is continuous at
$\overline{\alpha}$.

Theorem \ref{thm:2} is proved.\\

\textbf{Proof of theorem \ref{thm:3}}:\\

We define
$w_i=(-\triangle)^{i-1}u$, $w_{k+i}=(-\triangle)^{i-1}v, i=1,\ldots,k$.
Then $w$ satisfies:

\begin{equation} \label{eq:7}
\left\{ \begin{aligned}
         &-\triangle w_1 = w_2,  \cdots
             -\triangle w_{k-1}=w_k, \\
          &-\triangle w_{k} =w_{k+1}^p, \\
         &-\triangle w_{k+1}=w_{k+2},  \cdots
                -\triangle w_{2k-1}=w_{2k}, \\
          &-\triangle w_{2k}=w_{1}^q,  \\
          &w_1>0,\ldots,w_{2k}>0, \text{ in }{\mathbb{R}^n}
                          \end{aligned} \right.
                          \end{equation}
with $F(w)=(w_2,\ldots,w_k,w_1^p,w_{k+1},\ldots,w_{2k},w_{k+1}^q):\mathbb{R}_+^{2k}\rightarrow \mathbb{R}_+^{2k}$
continuous. $F(w)$ is also Lipschitz for $w>0$.
It is easy to check that $F$ is non-degenerate and $F(\beta) \neq 0$ when $\beta \neq 0$.
In \cite{LL13preprint}, Lei and Li have proved that the system:

\begin{equation} \label{eq:8}
\left\{ \begin{aligned}
          &-\triangle w_1 = w_2,  \cdots
             -\triangle w_{k-1}=w_k, \\
          &-\triangle w_{k} =w_{k+1}^p, \\
         &-\triangle w_{k+1}=w_{k+2},  \cdots
                -\triangle w_{2k-1}=w_{2k}, \\
          &-\triangle w_{2k}=w_{1}^q,  \\
          &w_1>0,\ldots,w_{2k}>0,  \text{ in } B_R(0)\\
          &w_1(x)=\ldots=w_{2k}(x)=0 \text{ on } \partial B_R(0)
                          \end{aligned} \right.
                          \end{equation}
admits no radially symmetric solution in the critical and super-critical cases
$\frac{1}{p+1}+\frac{1}{q+1}\leq\frac{n-2k}{n}$
for any $2k<n$.
Thus, according to Theorem \ref{thm:1}, system (\ref{eq:7})
admits a positive radial solution $w_i(x)>0$, $x \in\mathbb{R}^n$ , $i=1,\ldots,2k$.
Consequently, $u=w_1$ and $v=w_{k+1}$ solves system (\ref{eq:6}).

It is
interesting to point out that, for any solutions to (\ref{eq:6}), except the critical case where
$\frac{1}{p+1}+\frac{1}{q+1}=\frac{n-2k}{n}$, the `energy' is infinite in the sense that:
$\int_{\mathbb{R}^n} u^{q+1}=
\int_{\mathbb{R}^n} v^{p+1}=\int_{\mathbb{R}^n}
|(-\triangle)^{\frac{k}{2}} u|^2=\int_{\mathbb{R}^n} |(-\triangle)^{\frac{k}{2}} v|^2= \infty$.\\

 \textbf{Proof of theorem \ref{thm:4}}:\\

 The nonexistence in the subcritical case and the classification of solutions in the critical case
 have been proved in \cite{CLO06}. The existence follows exactly as the prove of theorem \ref{thm:3}.

 \paragraph{Acknowledgements.}
This work is partially supported by  NSFC grant 11271166 and by NSF grants DMS-0908097 and EAR-093467.

Congming Li

Department of Applied Mathematical, University of Colorado at Boulder,
Boulder, CO 80309, USA

Department of Mathematics, and MOE-LSC, Shanghai Jiao Tong University,
Shanghai, 200240, China


\begin{thebibliography}{20}

\bibitem{BrezisNirenberg83} H. Brezis and L. Nirenberg, 
    \emph{Positive solutions of nonlinear elliptic equations involving critical sobolev exponents},
    Comm. Pure Appl. Math.
    \textbf{36} (1983),  437--477.

\bibitem{CGS89}
    L. Caffarelli, B. Gidas, and J. Spruck,
    \emph{Asymptotic symmetry and local behavior of semilinear
    elliptic equations with critical Sobolev growth}, Comm. Pure
    Appl. Math. \textbf{42} (1989), 271--297.

\bibitem{CY97}
    A. Chang and P. Yang,
    \emph{On uniqueness of an n-th
    order differential equation in conformal geometry}, Math. Res.
    Letters,  \textbf{4} (1997), 1--12.

\bibitem{CL97}
    W. Chen and C. Li,
    \emph{A priori estimates for
    prescribing scalar curvature equations}, Ann. Math.
    \textbf{145} (1997), 547--564.

\bibitem{CL09}
    W. Chen and C. Li,
    \emph{An integral system and the Lane-Emden conjecture},
    Disc. Cont. Dynamics Sys.
    \textbf{24} (2009), 1167-1184.

\bibitem{CLO06}
    W. Chen, C. Li, and B. Ou,
    \emph{Classification
    of solutions for an integral equation}, Comm. Pure Appl.
    Math. \textbf{59} (2006), 330--343.

\bibitem{CEQ94} X. Chen, C. Elliott, and T. Qi,
    \emph{Shooting method for vortex solutions of a complex-valued Ginzburg-Landau equation},
    Proc. R. Soc. Edinb. 124(1994), 1075-1088.

\bibitem{Dancer01} E. N. Dancer,
    \emph{New solutions of equations on $\mathbb{R}^n$},
    Annali della Scuola Normale Superiore di Pisa - Classe di Scienze,
    S¨¦r. 4, 30 no. 3-4 (2001), 535-563.

\bibitem{Gui96} C. Gui, \emph{On positive entire solutions of the elliptic equation $¦¤u+K(x)u^p=0$ and
        its applications to Riemannian geometry}, Proc. Roy. Soc. Edinburgh Sect. A, 126 (1996),  225--237.

\bibitem{JL73} D. Joseph and T. Lundgren,
    \emph{Quasilinear Dirichlet problem driven by positive sources},
    Arch. Rational Mech. Anal. 49 (1972) 241-269.

\bibitem{LL13preprint} Y. Lei and C. Li,
    \emph{Sharp criteria of Liouville type for some nonlinear systems},
    arXiv:1301.6235.

\bibitem{Lieb83}
    E. Lieb,
    \emph{Sharp constants in the Hardy-Littlewood-Sobolev and
    related inequalities}, Ann. Math. \textbf{118} (1983),
    349--374.

\bibitem{NS86}
    W.-M. Ni and J. Serrin,
    \emph{Existence and nonexistence theorems for ground states of quasilinear partial differential
    equations. The anomalous case}, Accad. Naz. Lincei. \textbf{77} (1986), 231-257.

\bibitem{Nirenberg01} L. Nirenberg, \emph{Topics in Nonlinear Functional Analysis}. Notes by R. A. Artino,
    Vol. 6 of Courant Lecture Notes in Mathematics, New York: New York University Courant Institute of Mathematical Sciences 2001.

\bibitem{PS86}
    P. Pucci and J. Serrin,
    \emph{A general variational identity},
    Indiana Univ. J. Math. \textbf{35} (1986), 681-703.

\bibitem{PV08}
      N. Phuc and I. Verbitsky,
      \emph{Quasilinear and Hessian equations of Lane-Emden type},
      Ann.  Math. \textbf{168}
    (2008), 859--914.


\bibitem{Schwartz69} J. T, Schwartz, \emph{Nonlinear Functional Analysis}. Notes by H. Fattorni,
    R. Nirenberg, and H. Porta, with an additional chapter by H. Karcher. Notes on Mathematics and Its
    Applications. Gordon and Breach, New York-London-Paris, 1969.

\bibitem{SZ96}
     J. Serrin and H. Zou,
     \emph{Non-existence of positive solution of Lane-Emden systems},
     Differential and Integral equations, \textbf{9} (1996), 635--653.

\bibitem{SZ02}
     J. Serrin and H. Zou,
     \emph{Cauchy-Liouville and universal boundedness theorems for
     quasilinear elliptic equations and inequalities},
    Acta Math.\textbf{189} (2002), 79--142.

\bibitem{Souplet09}
     P. Souplet,
     \emph{The proof of the Lane-Emden conjecture in 4 space dimensions},
     Adv. Math. \textbf{221} (2009), 1409--1427.

\end{thebibliography}
\end{document}